\documentclass[12pt]{article}
\usepackage{amsmath,amssymb,amsthm}
\usepackage{vmargin}
\usepackage{url}
\usepackage{nicefrac}

\usepackage{graphicx}

\usepackage{algpseudocode}

\usepackage{rotating}  

\setpapersize{USletter}
\setmargnohfrb{1.0in}{1.0in}{1.0in}{1.0in}

\AtBeginDocument{
\footskip=30pt
}

\usepackage{pxfonts}

\newcommand{\beq}{\begin{equation}}
\newcommand{\eeq}{\end{equation}}
\newcommand{\bea}{\begin{array}}
\newcommand{\eea}{\end{array}}

\newcommand{\B}{\mathcal{B}}
\newcommand{\U}{\mathcal{U}}

\title{On the number of permutations with bounded run lengths}
\author{Max A. Alekseyev\thanks{Department of Computer Science and Engineering,
University of South Carolina, Columbia, SC, U.S.A.\newline{Email: maxal@cse.sc.edu}} }

\begin{document}
\maketitle

\begin{abstract}
In this work we obtain recurrent formulae for the number of permutations with 
either increasing or monotonic (i.e., both increasing and decreasing) 
runs of bounded length.
Our formulae allow one to efficiently compute the number of such permutations.
In particular, we use the formulae to find and correct a few miscalculations in the classic 1966 book by David, Kendall, and Barton.

We further use our formulae to derive differential equations for the corresponding exponential generating functions.
In the case of increasing runs, we solve these equations and obtain closed-form expressions for the generating functions.
\end{abstract}

\section{Introduction}

A (monotonic) \emph{run} in a permutation $p=(p_1,p_2,\dots,p_n)$ is a maximal increasing or decreasing subsequence of consecutive elements in $p$.
Similarly, an \emph{increasing} (resp. \emph{decreasing}) \emph{run} in $p$ is a maximal increasing (resp. decreasing) subsequence of consecutive elements in $p$.
David et al.~\cite{David1966} in Tables 7.4.1 and 7.4.2 give counts\footnote{Colin Mallows and Neil Sloane suggest
that these counts were almost certainly the result of hand calculations.} for
\begin{itemize}
\item the number $I^k(n)$ of order $n$ permutations whose longest increasing run length equals $k$ (for $k,n\leq 18$);
\item the number $A^k(n)$ of order $n$ permutations whose longest monotonic run length equals $k$ (for $k,n\leq 14$).
\end{itemize}
It turns out that their counts for $I^k(n)$ and $A^k(n)$ are incorrect for $n\geq 16$ and $n\geq 13$ respectively.

We notice that $I^k(n) = U^k(n) - U^{k-1}(n)$ and $A^k(n) = B^k(n) - B^{k-1}(n)$ where $U^k(n)$ and $B^k(n)$ 
is the number of permutations of order $n$ whose runs length does not exceed $k$.

In this note we derive recurrent formulae for $U^k(n)$ and $B^k(n)$ 
as well as differential equations for their exponential generating functions.
These formulae allowed us to compute $B^k(n)$ and $A^k(n)$ accurately and correct miscalculations in~\cite{David1966}.

We remark that another way to obtain differential equations 
for these generating functions was described by Elizalde and Noy~\cite{Elizalde2003} who
studied a more general problem of counting permutations with forbidden subpermutations
using symbolic methods. In contrast, we use only elementary observations to
obtain recurrent formulae for $U^k(n)$ and $B^k(n)$ and describe algorithms for computing them. 
We also derive an explicit closed-form expressions for the generating functions for $U^k(n)$.

\section{Recurrent formulae for $U^k(n)$ and $B^k(n)$}

Let $p=(p_1,p_2,\dots,p_n)$ be a permutation of order $n>1$ with runs length not exceeding $k$. 
It is easy to see that $p_t$ may be a maximum element of $p$ (i.e., $p_t=n$) only in one of 
the following three cases:
\begin{enumerate}
\item $t=1$ and $p_1 > p_2$; 
\item $t=n$ and $p_{n-1} < p_n$; 
\item $1< t <n$, $p_{t-1} < p_t$, and $p_t > p_{t+1}$. 
\end{enumerate}

In case 1, we remove the first element from $p$ to obtain a permutation 
$p'=(p_2,p_3,\dots,p_n)$ of order $n-1$.
The permutation $p'$ here can be any permutation whose (monotonic or increasing) runs length does not exceed $k$
with an additional restriction (in the case of bounded monotonic runs) that 
if the initial run in $p'$ is decreasing, 
then its length does not exceed $k-1$.

Similarly, in case 2, we remove the last element from $p$ to obtain a permutation
$p'=(p_1,p_2,\dots,p_{n-1})$ of order $n-1$.
The permutation $p'$ here can be any permutation whose (monotonic or increasing) runs length does not exceed $k$
with an additional restriction that
if the final run in $p'$ is increasing 
then its length does not exceed $k-1$.

In case 3, removing of the element $p_t$ splits $p$ into two vectors
$(p_1,p_2,\dots,p_{t-1})$ and $(p_{t+1},p_{t+2},\dots,p_n)$. 
We relabel their elements with integers $1,2,\dots,t-1$ and $1,2,\dots,n-t$ 
(preserving the order relationship) to obtain permutations $p'$ and $p"$ of order $t-1$ and $n-t$ respectively.
By construction, the permutation $p'$ has the same length and order of runs as the prefix of length $t-1$ of $p$, 
while the permutations $p"$ has the same length and order of runs as the suffix of length $n-t$ of $p$.
Therefore, $p'$ and $p"$ can be any permutations whose runs length does not exceed $k$ with
additional restrictions that 
\begin{itemize}
\item if the final run in $p'$ in increasing, then its length is at most $k-1$;
\item (in case of bounded monotonic runs) if the initial run in $p"$ is decreasing, then its length is at most $k-1$.
\end{itemize}

Let $U^k_j$ be the number of permutations $p$ of order $n$ 
whose increasing runs length does not exceed $k$, and
the final increasing run (if it is present) in $p$ has length at most $j$.
Trivially $U^k_j(1) = 1$ for any $1\leq j\leq k$. 
We find it convenient to define $U^k_j(n) = 0$ whenever $j<1$.

Similarly, let $B^k_{i,j}(n)$ be the number of permutations $p$ of order $n$ 
whose runs length does not exceed $k$, and
the initial decreasing run (if it is present) in $p$  has length at most $i$ and the final increasing run (if it is present) in $p$ has length at most $j$.
Trivially $B^k_{i,j}(1) = 1$ for any $1\leq i,j\leq k$. 
We find it convenient to define $B^k_{i,j}(n) = 0$ whenever $i<1$ or $j<1$.

The above observations lead to the following formulae:
$$U^k_j(n) = U^k_j(n-1) + U^k_{j-1}(n-1) + \sum_{t=2}^{n-1} \binom{n-1}{t-1}\cdot U^k_{k-1}(t-1)\cdot U^k_j(n-t)$$
$$B^k_{i,j}(n) = B^k_{i-1,j}(n-1) + B^k_{i,j-1}(n-1) + \sum_{t=2}^{n-1} \binom{n-1}{t-1}\cdot B^k_{i,k-1}(t-1)\cdot B^k_{k-1,j}(n-t)$$
which hold for $n>1$ and any $1\leq i,j\leq k$.
Here the binomial coefficient $\tbinom{n-1}{t-1}$ stands for the number 
of ways to distribute elements $1,2,\dots,n-1$ of $p$ between the prefix and 
suffix corresponding to the permutations $p'$ and $p"$ in case 3.

We also remark that the involution $(p_1,p_2,\dots,p_n)\mapsto (p_n,p_{n-1},\dots,p_1)$ on the set of all permutations of order $n$
implies that $B^k_{i,j}(n) = B^k_{j,i}(n)$.

\section{Computing $U^k(n)$ and $B^k(n)$}

Suppose that $k$ is fixed. It is easy to see that 
$U^k(n) = U^k_k(n)$ and $B^k(n) = B^k_{k,k}(n)$ which allows one
to compute them efficiently using the recurrent formulae for $U^k_j(n)$ and $B^k_{i,j}(n)$.
In particular, $U^k(n) = U^k_k(n)$ for all $n\leq N$ can be computed as follows:
\begin{algorithmic}[1]
\State $U\gets$ an array of size $k$
\For{$j\gets 1$ to $k$}
\State $U[j] \gets$ an array of integers of size $N$
\State $U[j][1] \gets 1$
\EndFor
\For{$n\gets 2$ to $N$}
 \For{$j\gets 1$ to $k$}
        \State $U[j][n] \gets U[j][n-1]$
        \If{$j>1$} 
        \State $U[j][n] \gets U[j][n] + U[j-1][n-1]$ 
        \EndIf
        \If{$n>2$ and $k>1$}
        \State $U[j][n] \gets U[j][n] + \sum_{t=2}^{n-1} \tbinom{n-1}{t-1}\cdot U[k-1][t-1]\cdot U[j][n-t]$ 
        \EndIf
      \EndFor
\EndFor
\State \Return $U[k]$
\end{algorithmic}

Similarly, $B^k(n) = B^k_{k,k}(n)$ for all $n\leq N$ can be computed as follows:
\begin{algorithmic}[1]
\State $B\gets$ an array of size $k\times k$
\For{$i\gets 1$ to $k$}
\For{$j\gets 1$ to $k$}
\State $B[i,j]\gets$ an array of integers of size $N$
\State $B[i,j][1] \gets 1$
\EndFor
\EndFor
\For{$n\gets 2$ to $N$}
  \For{$i\gets 1$ to $k$}
    \For{$j\gets 1$ to $k$}
      \State $B[i,j][n] \gets 0$
      \If{$i>1$}
        \State $B[i,j][n] \gets B[i,j][n] + B[i-1,j][n-1]$
      \EndIf
      \If{$j>1$}
        \State $B[i,j][n] \gets B[i,j][n] + B[i,j-1][n-1]$
      \EndIf
      \If{$n>2$ and $k>1$}
        \State $B[i,j][n] \gets B[i,j][n] + \sum_{t=2}^{n-1} \tbinom{n-1}{t-1}\cdot B[i,k-1][t-1]\cdot B[k-1,j][n-t]$
      \EndIf
    \EndFor
  \EndFor
\EndFor
\State \Return $B[k,k]$.
\end{algorithmic}



We used these algorithms to compute values $U^k(n)$ and $B^k(n)$ for $k,n\leq 18$ and listed them in Tables~\ref{table:Ukn} and \ref{table:Bkn} respectively.
Subtracting from each row the previous one, we obtain Tables~\ref{table:Ikn} and \ref{table:Akn} listing values of $I^k(n)$ and $A^k(n)$.
We remark that Table~\ref{table:Ikn} is present (column-wise) in the OEIS~\cite{OEIS} as sequence \textsc{A008304} with 
its rows (for $2\leq k\leq 6$) given by sequences \textsc{A008303}, \textsc{A000402}, \textsc{A000434}, \textsc{A000456}, and \textsc{A000467}.
Table~\ref{table:Akn} (column-wise) is present in the OEIS as sequence \textsc{A211318} with its rows (for $2\leq k\leq 5$) given by sequences \textsc{A001250}, \textsc{A001251}, 
\textsc{A001252}, and \textsc{A001253}.



\section{Exponential generating function for $U^k(n)$}

For fixed integers $k,j$, let $\U_j^k(x)$ be the exponential generating function for $U^k_j(n)$:
$$\U_j^k(x) = \sum_{n=0}^{\infty} U^k_j(n)\cdot \frac{x^n}{n!}.$$

The recurrent formula for $U^k_j(n)$ implies the following system of differential equations:
$$\left\{\frac{d}{dx}\U_j^k(x) = 1 + \U_j^k(x) + \U_{j-1}^k(x) + \U_{k-1}^k(x)\cdot \U_{j}^k(x),\right.
\qquad j=1,2,\dots,k.$$ 

In particular, for $j=k$ we have
$$\frac{d}{dx}\U_k^k(x) = (1 + \U_k^k(x))\cdot (1 + \U_{k-1}^k(x))$$
or
$$1 + \U_{k-1}^k(x) = \frac{1}{1 + \U_k^k(x)} \frac{d}{dx}\U_k^k(x).$$
Plugging this into the $j$-th equation, we conclude 
$$\frac{d}{dx}\left(\frac{\U_j^k(x)}{1 + \U_k^k(x)}\right) 
= \frac{1}{1+\U_k^k(x)} + \frac{\U_{j-1}^k(x)}{1+\U_k^k(x)}$$
implying that for $y(x) = \tfrac{1}{1 + \U_k^k(x)}$ and any $j=0,2,\dots,k-1$, we have
$$(1-y)^{(k-j)} = y^{(k-1-j)} + y^{(k-2-j)} + \cdots + y' + y + y\cdot \U_j^k(x).$$

In particular, for $j=0$ we get the following linear differential equation:
$$y^{(k)} + y^{(k-1)} + y^{(k-2)} + \cdots + y' + y = 0.$$
General solution to this equation is
$$y(x) = \sum_{i=1}^{k} c_i\cdot e^{r^i\cdot x}$$
where $r$ is a primitive $(k+1)$-st degree root of $1$ and $c_i$ are constant coefficients
such that $\sum_{i=1}^k c_i = y(0) = 1$.

For $j=1,2,\dots,k-1$, we expressing $\U_j^k(x)$ in terms of $y$ as follows:
$$\U_j^k(x) = - \frac{1}{y}\cdot \sum_{m=0}^{k-j} y^{(m)}.$$
Since $\U_j^k(0) = 0$, this expression implies
$$\sum_{m=0}^{k-j} \sum_{i=1}^{k} c_i\cdot r^{im} = 0
\qquad\text{or}\qquad
\sum_{i=1}^{k} c_i\cdot \frac{r^{i(k-j+1)}-1}{r^i - 1} = 0.$$
It is easy to see that solution to this system of linear equations is $c_i = \frac{r^i - 1}{(k+1)r^i}$, implying that
$$
y(x) = \sum_{n=0}^{\infty} q_n \frac{x^n}{n!}\qquad
\text{and}\qquad
\U_k^k(x) = \left( \sum_{n=0}^{\infty} q_n \frac{x^n}{n!} \right)^{-1} - 1
$$
where
$$
q_n = \begin{cases} 
1 & \text{if}\ n\equiv 0\pmod{k+1},\\
-1 & \text{if}\ n\equiv 1\pmod{k+1},\\
0 & \text{if}\ n\not{\equiv} 0,1\pmod{k+1}.
\end{cases}
$$


\subsection{Examples}

Exponential generating functions of $U^k(n)$ for $k=1,2,3,4$ are:
$$\U_1^1(x) = \left( 1 - \frac{x}{1!} + \frac{x^2}{2!} - \frac{x^3}{3!} + \ldots\right)^{-1} - 1 =
1\cdot \frac{x}{1} + 1\cdot \frac{x^2}{2!} + 1\cdot \frac{x^3}{3!} + 1\cdot \frac{x^4}{4!} + 1\cdot\frac{x^5}{5!} +\ldots$$

$$\U_2^2(x) = \left( 1 - \frac{x}{1!} + \frac{x^3}{3!} - \frac{x^4}{4!} + \ldots\right)^{-1} - 1 =
1\cdot\frac{x}{1} + 2\cdot \frac{x^2}{2!} + 5\cdot \frac{x^3}{3!} + 17\cdot \frac{x^4}{4!} + 70\cdot\frac{x^5}{5!} +\ldots$$

$$\U_3^3(x) = \left( 1 - \frac{x}{1!} + \frac{x^4}{4!} - \frac{x^5}{5!} + \ldots\right)^{-1} - 1 =
1\cdot\frac{x}{1} + 2\cdot \frac{x^2}{2!} + 6\cdot \frac{x^3}{3!} + 23\cdot \frac{x^4}{4!} + 111\cdot\frac{x^5}{5!} +\ldots$$

$$\U_4^4(x) = \left( 1 - \frac{x}{1!} + \frac{x^5}{5!} - \frac{x^6}{6!} + \ldots\right)^{-1} - 1 =
1\cdot\frac{x}{1} + 2\cdot \frac{x^2}{2!} + 6\cdot \frac{x^3}{3!} + 24\cdot \frac{x^4}{4!} + 119\cdot\frac{x^5}{5!} +\ldots$$

The difference $\U_k^k(x)-\U_{k-1}^{k-1}(x)$ gives the exponential generating function of the numbers $I^k(n)$:
$$\U_2^2(x) - \U_1^1(x) = 
0\cdot \frac{x}{1} + 1\cdot \frac{x^2}{2!} + 4\cdot \frac{x^3}{3!} + 16\cdot \frac{x^4}{4!} + 69\cdot\frac{x^5}{5!} +\ldots$$

$$\U_3^3(x) - \U_2^2(x) = 
0\cdot \frac{x}{1} + 0\cdot \frac{x^2}{2!} + 1\cdot \frac{x^3}{3!} + 6\cdot \frac{x^4}{4!} + 41\cdot\frac{x^5}{5!} +\ldots$$

$$\U_4^4(x) - \U_3^3(x) = 
0\cdot \frac{x}{1} + 0\cdot \frac{x^2}{2!} + 0\cdot \frac{x^3}{3!} + 1\cdot \frac{x^4}{4!} + 8\cdot\frac{x^5}{5!} +\ldots$$


\section{Exponential generating function for $B^k(n)$}

For fixed integers $k,i,j$, let $\B_{i,j}^k(x)$ be the exponential generating function for $B^k_{i,j}(n)$:
$$\B_{i,j}^k(x) = \sum_{n=0}^{\infty} B^k_{i,j}(n)\cdot \frac{x^n}{n!}.$$

The recurrent formula for $B^k_{i,j}(n)$ implies the following system of differential equations:
$$\left\{\frac{d}{dx}\B_{i,j}^k(x) = 1 + \B_{i-1,j}^k(x) + \B_{i,j-1}^k(x) + \B_{i,k-1}^k(x)\cdot \B_{k-1,j}^k(x),\right.
\qquad i,j=1,2,\dots,k.$$ 
Because of the symmetry $\B_{i,j}^k(x) = \B_{j,i}^k(x)$,
this system contains $\frac{k(k+1)}{2}$ distinct functions and equations.

\subsection{Examples}

\subsubsection{$k=2$}

For $k=2$, we have the following system of differential equations:
$$
\begin{cases}
\frac{d}{dx}\B_{2,2}^2(x) = 1 + 2\B_{1,2}^2(x) + \B_{1,2}^2(x)^2, \\
\frac{d}{dx}\B_{1,2}^2(x) = 1 + \B_{1,1}^2(x) + \B_{1,1}^2(x)\cdot \B_{1,2}^2(x), \\
\frac{d}{dx}\B_{1,1}^2(x) = 1 + \B_{1,1}^2(x)^2.
\end{cases}
$$
with the constraints $\B_{i,j}^k(0) = 0$.

The system has the following solution:
$$
\begin{cases}
\B_{1,1}^2(x) = \tan(x), \\
\B_{1,2}^2(x) = \tan(x) + \sec(x) - 1,\\
\B_{2,2}^2(x) =  2(\tan(x) + \sec(x)-1) - x.
\end{cases}
$$
Therefore, the numbers $B^2(n)$ represent the coefficients in the series expansion
$$\B_{2,2}^2(x) = 2(\tan(x) + \sec(x)-1) - x = x + 2\cdot \frac{x^2}{2!} + 4\cdot \frac{x^3}{3!} + 10\cdot \frac{x^4}{4!} + 32\cdot \frac{x^5}{5!} + \cdots$$
and form sequence \textsc{A001250} in~\cite{OEIS}. It also counts the number of permutations of order $n$ with exactly $n-1$ runs~\cite{Comtet1974}.

\subsubsection{$k=3$}

For $k=3$, we have the following system of differential equations:
$$
\begin{cases}
\frac{d}{dx}\B_{3,3}^3(x) = 1 + 2\B_{2,3}^3(x) + \B_{2,3}^3(x)^2 \\
\frac{d}{dx}\B_{2,3}^3(x) = 1 + \B_{1,3}^3(x) + \B_{2,2}^3(x) + \B_{2,2}^3(x)\cdot \B_{2,3}^3(x) \\
\frac{d}{dx}\B_{1,3}^3(x) = 1 + \B_{1,2}^3(x) + \B_{1,2}^3(x)\cdot \B_{2,3}^3(x) \\
\frac{d}{dx}\B_{2,2}^3(x) = 1 + 2\B_{1,2}^3(x) + \B_{2,2}^3(x)^2 \\
\frac{d}{dx}\B_{1,2}^3(x) = 1 + \B_{1,1}^3(x) + \B_{1,2}^3(x)\cdot \B_{2,2}^3(x) \\
\frac{d}{dx}\B_{1,1}^3(x) = 1 + \B_{1,2}^3(x)^2 
\end{cases}
$$

Although we have not been able to solve this system, we remark that the last three equations form a subsystem involving the 
functions $\B_{2,2}^3(x)$, $\B_{1,2}^3(x)$, and $\B_{1,1}^3(x)$,
which implies the following autonomous ordinary differential equation for $y(x) = \B_{2,2}^3(x)$:
$$2\cdot y"' - 6\cdot y\cdot y" - 7\cdot y'^2 + 8\cdot y^2\cdot y' + 4\cdot y' - y^4 - 2\cdot y^2 - 5 = 0.$$
This equation further reduces to the following second order differential equation for $w(y) = y'$ \cite{Polyanin2003}:
$$2\cdot w^2\cdot w" + 2\cdot w\cdot w'^2 - 6\cdot y\cdot w\cdot w' 
- 7\cdot w^2 + 8 \cdot y^2 \cdot w + 4\cdot w - y^4 - 2\cdot y^2 - 5 = 0.$$
Solving this equation would be a step towards obtaining the generating function $\B^3_{3,3}(x)$ for the numbers $B^3(n)$.

\section*{Acknowledgements}

The author thanks Sean A. Irvine for raising concerns about correctness of the values $A^k(n)$ listed in \cite{David1966} 
and posing the problem of computing $A^k(n)$ efficiently. The author is also thankful to Neil Sloane for a number of helpful comments.

\begin{footnotesize}
\bibliographystyle{plain} 
\bibliography{perm.bib} 
\end{footnotesize}

\newpage
\clearpage

\begin{sidewaystable}
\begin{center}
\begin{tiny}
\caption{Values $U^k(n)$ for $k,n\leq 18$.}\label{table:Ukn}
\begin{tabular}{|r|rrrrrrrrrrrrrrrrrr|}
\hline
$n$ & 1 & 2 & 3 & 4 & 5 & 6 & 7 & 8 & 9 & 10 & 11 & 12 & 13 & 14 & 15 & 16 & 17 & 18 \\
\hline\hline
$U^{1}(n)$ & 1 & 1 & 1 & 1 & 1 & 1 & 1 & 1 & 1 & 1 & 1 & 1 & 1 & 1 & 1 & 1 & 1 & 1 \\
\hline
$U^{2}(n)$ & 1 & 2 & 5 & 17 & 70 & 349 & 2017 & 13358 & 99377 & 822041 & 7477162 & 74207209 & 797771521 & 9236662346 & 114579019469 & 1516103040833 & 21314681315998 & 317288088082405 \\
\hline
$U^{3}(n)$ & 1 & 2 & 6 & 23 & 111 & 642 & 4326 & 33333 & 288901 & 2782082 & 29471046 & 340568843 & 4263603891 & 57482264322 & 830335952166 & 12793889924553 & 209449977967081 & 3630626729775362 \\
\hline
$U^{4}(n)$ & 1 & 2 & 6 & 24 & 119 & 709 & 4928 & 39144 & 349776 & 3472811 & 37928331 & 451891992 & 5832672456 & 81074690424 & 1207441809209 & 19181203110129 & 323753459184738 & 5785975294622694 \\
\hline
$U^{5}(n)$ & 1 & 2 & 6 & 24 & 120 & 719 & 5027 & 40168 & 361080 & 3606480 & 39623760 & 474915803 & 6166512899 & 86227808578 & 1291868401830 & 20645144452320 & 350547210173280 & 6302294420371031 \\
\hline
$U^{6}(n)$ & 1 & 2 & 6 & 24 & 120 & 720 & 5039 & 40305 & 362682 & 3626190 & 39881160 & 478490760 & 6219298800 & 87055051511 & 1305598835941 & 20885951018102 & 354999461960226 & 6388879812001704 \\
\hline
$U^{7}(n)$ & 1 & 2 & 6 & 24 & 120 & 720 & 5040 & 40319 & 362863 & 3628550 & 39913170 & 478947480 & 6226179960 & 87164597520 & 1307440134000 & 20918580896069 & 355608034188517 & 6400803479701178 \\
\hline
$U^{8}(n)$ & 1 & 2 & 6 & 24 & 120 & 720 & 5040 & 40320 & 362879 & 3628781 & 39916492 & 478996716 & 6226941864 & 87176969880 & 1307651304960 & 20922368987520 & 355679390626560 & 6402213152423659 \\
\hline
$U^{9}(n)$ & 1 & 2 & 6 & 24 & 120 & 720 & 5040 & 40320 & 362880 & 3628799 & 39916779 & 479001228 & 6227014404 & 87178179816 & 1307672369640 & 20922752672640 & 355686706327680 & 6402359109968640 \\
\hline
$U^{10}(n)$ & 1 & 2 & 6 & 24 & 120 & 720 & 5040 & 40320 & 362880 & 3628800 & 39916799 & 479001577 & 6227020358 & 87178283010 & 1307674215120 & 20922786961440 & 355687370176320 & 6402372516146880 \\
\hline
$U^{11}(n)$ & 1 & 2 & 6 & 24 & 120 & 720 & 5040 & 40320 & 362880 & 3628800 & 39916800 & 479001599 & 6227020775 & 87178290682 & 1307674357710 & 20922789683040 & 355687423926240 & 6402373618334400 \\
\hline
$U^{12}(n)$ & 1 & 2 & 6 & 24 & 120 & 720 & 5040 & 40320 & 362880 & 3628800 & 39916800 & 479001600 & 6227020799 & 87178291173 & 1307674367400 & 20922789875280 & 355687427826720 & 6402373699926240 \\
\hline
$U^{13}(n)$ & 1 & 2 & 6 & 24 & 120 & 720 & 5040 & 40320 & 362880 & 3628800 & 39916800 & 479001600 & 6227020800 & 87178291199 & 1307674367971 & 20922789887312 & 355687428080496 & 6402373705380384 \\
\hline
$U^{14}(n)$ & 1 & 2 & 6 & 24 & 120 & 720 & 5040 & 40320 & 362880 & 3628800 & 39916800 & 479001600 & 6227020800 & 87178291200 & 1307674367999 & 20922789887969 & 355687428095218 & 6402373705709334 \\
\hline
$U^{15}(n)$ & 1 & 2 & 6 & 24 & 120 & 720 & 5040 & 40320 & 362880 & 3628800 & 39916800 & 479001600 & 6227020800 & 87178291200 & 1307674368000 & 20922789887999 & 355687428095967 & 6402373705727118 \\
\hline
$U^{16}(n)$ & 1 & 2 & 6 & 24 & 120 & 720 & 5040 & 40320 & 362880 & 3628800 & 39916800 & 479001600 & 6227020800 & 87178291200 & 1307674368000 & 20922789888000 & 355687428095999 & 6402373705727965 \\
\hline
$U^{17}(n)$ & 1 & 2 & 6 & 24 & 120 & 720 & 5040 & 40320 & 362880 & 3628800 & 39916800 & 479001600 & 6227020800 & 87178291200 & 1307674368000 & 20922789888000 & 355687428096000 & 6402373705727999 \\
\hline
$U^{18}(n)$ & 1 & 2 & 6 & 24 & 120 & 720 & 5040 & 40320 & 362880 & 3628800 & 39916800 & 479001600 & 6227020800 & 87178291200 & 1307674368000 & 20922789888000 & 355687428096000 & 6402373705728000 \\
\hline
\end{tabular}

\vspace{2\baselineskip}
\caption{Values $I^k(n)$ for $k,n\leq 18$. Highlighted values indicate disagreements with \cite{David1966}.}\label{table:Ikn}
\begin{tabular}{|r|rrrrrrrrrrrrrrrrrr|}
\hline
$n$ & 1 & 2 & 3 & 4 & 5 & 6 & 7 & 8 & 9 & 10 & 11 & 12 & 13 & 14 & 15 & 16 & 17 & 18 \\
\hline\hline
$I^1(n)$ & 1 & 1 & 1 & 1 & 1 & 1 & 1 & 1 & 1 & 1 & 1 & 1 & 1 & 1 & 1 & 1 & 1 & 1 \\
\hline
$I^2(n)$ & 0 & 1 & 4 & 16 & 69 & 348 & 2016 & 13357 & 99376 & 822040 & 7477161 & 74207208 & 797771520 & 9236662345 & 114579019468 & 1516103040832 & 21314681315997 & 317288088082404 \\
\hline
$I^3(n)$ & 0 & 0 & 1 & 6 & 41 & 293 & 2309 & 19975 & 189524 & 1960041 & 21993884 & 266361634 & 3465832370 & 48245601976 & 715756932697 & \textbf{11277786883720} & \textbf{188135296651083} & \textbf{3313338641692957} \\
\hline
$I^4(n)$ & 0 & 0 & 0 & 1 & 8 & 67 & 602 & 5811 & 60875 & 690729 & 8457285 & 111323149 & 1569068565 & 23592426102 & 377105857043 & \textbf{6387313185576} & \textbf{114303481217657} & \textbf{2155348564847332} \\
\hline
$I^5(n)$ & 0 & 0 & 0 & 0 & 1 & 10 & 99 & 1024 & 11304 & 133669 & 1695429 & 23023811 & 333840443 & 5153118154 & 84426592621 & 1463941342191 & 26793750988542 & 516319125748337 \\
\hline
$I^6(n)$ & 0 & 0 & 0 & 0 & 0 & 1 & 12 & 137 & 1602 & 19710 & 257400 & 3574957 & 52785901 & 827242933 & 13730434111 & 240806565782 & 4452251786946 & 86585391630673 \\
\hline
$I^7(n)$ & 0 & 0 & 0 & 0 & 0 & 0 & 1 & 14 & 181 & 2360 & 32010 & 456720 & 6881160 & 109546009 & 1841298059 & 32629877967 & 608572228291 & 11923667699474 \\
\hline
$I^8(n)$ & 0 & 0 & 0 & 0 & 0 & 0 & 0 & 1 & 16 & 231 & 3322 & 49236 & 761904 & 12372360 & 211170960 & 3788091451 & 71356438043 & 1409672722481 \\
\hline
$I^9(n)$ & 0 & 0 & 0 & 0 & 0 & 0 & 0 & 0 & 1 & 18 & 287 & 4512 & 72540 & 1209936 & 21064680 & 383685120 & 7315701120 & 145957544981 \\
\hline
$I^{10}(n)$ & 0 & 0 & 0 & 0 & 0 & 0 & 0 & 0 & 0 & 1 & 20 & 349 & 5954 & 103194 & 1845480 & 34288800 & 663848640 & 13406178240 \\
\hline
$I^{11}(n)$ & 0 & 0 & 0 & 0 & 0 & 0 & 0 & 0 & 0 & 0 & 1 & 22 & 417 & 7672 & 142590 & 2721600 & 53749920 & 1102187520 \\
\hline
$I^{12}(n)$ & 0 & 0 & 0 & 0 & 0 & 0 & 0 & 0 & 0 & 0 & 0 & 1 & 24 & 491 & 9690 & 192240 & 3900480 & 81591840 \\
\hline
$I^{13}(n)$ & 0 & 0 & 0 & 0 & 0 & 0 & 0 & 0 & 0 & 0 & 0 & 0 & 1 & 26 & 571 & 12032 & 253776 & 5454144 \\
\hline
$I^{14}(n)$ & 0 & 0 & 0 & 0 & 0 & 0 & 0 & 0 & 0 & 0 & 0 & 0 & 0 & 1 & 28 & 657 & 14722 & 328950 \\
\hline
$I^{15}(n)$ & 0 & 0 & 0 & 0 & 0 & 0 & 0 & 0 & 0 & 0 & 0 & 0 & 0 & 0 & 1 & 30 & 749 & 17784 \\
\hline
$I^{16}(n)$ & 0 & 0 & 0 & 0 & 0 & 0 & 0 & 0 & 0 & 0 & 0 & 0 & 0 & 0 & 0 & 1 & 32 & 847 \\
\hline
$I^{17}(n)$ & 0 & 0 & 0 & 0 & 0 & 0 & 0 & 0 & 0 & 0 & 0 & 0 & 0 & 0 & 0 & 0 & 1 & 34 \\
\hline
$I^{18}(n)$ & 0 & 0 & 0 & 0 & 0 & 0 & 0 & 0 & 0 & 0 & 0 & 0 & 0 & 0 & 0 & 0 & 0 & 1 \\
\hline
\end{tabular}
\end{tiny}
\end{center}
\end{sidewaystable}

\begin{sidewaystable}
\begin{center}
\begin{tiny}
\caption{Values $B^k(n)$ for $k,n\leq 18$.}\label{table:Bkn}
\begin{tabular}{|r|rrrrrrrrrrrrrrrrrr|}
\hline
$n$ & 1 & 2 & 3 & 4 & 5 & 6 & 7 & 8 & 9 & 10 & 11 & 12 & 13 & 14 & 15 & 16 & 17 & 18 \\
\hline\hline
$B^{1}(n)$ & 1 & 0 & 0 & 0 & 0 & 0 & 0 & 0 & 0 & 0 & 0 & 0 & 0 & 0 & 0 & 0 & 0 & 0 \\
\hline
$B^{2}(n)$ & 1 & 2 & 4 & 10 & 32 & 122 & 544 & 2770 & 15872 & 101042 & 707584 & 5405530 & 44736512 & 398721962 & 3807514624 & 38783024290 & 419730685952 & 4809759350882 \\
\hline
$B^{3}(n)$ & 1 & 2 & 6 & 22 & 102 & 564 & 3652 & 26986 & 224458 & 2073946 & 21080922 & 233752052 & 2807949492 & 36324988206 & 503484183878 & 7443797211854 & 116931715588046 & 1944883690208684 \\
\hline
$B^{4}(n)$ & 1 & 2 & 6 & 24 & 118 & 698 & 4816 & 37968 & 336812 & 3319622 & 35990262 & 425668584 & 5454050314 & 75257838602 & 1112621686120 & 17545752570360 & 293985178842320 & 5215578061637498 \\
\hline
$B^{5}(n)$ & 1 & 2 & 6 & 24 & 120 & 718 & 5014 & 40016 & 359280 & 3584160 & 39331224 & 470842102 & 6106259878 & 85282508228 & 1276168085580 & 20369694217750 & 345453884789910 & 6203249305454148 \\
\hline
$B^{6}(n)$ & 1 & 2 & 6 & 24 & 120 & 720 & 5038 & 40290 & 362484 & 3623580 & 39845520 & 477979920 & 6211578648 & 86931863566 & 1303524552206 & 20849140937272 & 354312156550056 & 6375401280887904 \\
\hline
$B^{7}(n)$ & 1 & 2 & 6 & 24 & 120 & 720 & 5040 & 40318 & 362846 & 3628300 & 39909540 & 478893360 & 6225339120 & 87150903840 & 1307205906864 & 20914372123786 & 355528646248138 & 6399233407501172 \\
\hline
$B^{8}(n)$ & 1 & 2 & 6 & 24 & 120 & 720 & 5040 & 40320 & 362878 & 3628762 & 39916184 & 478991832 & 6226862928 & 87175648560 & 1307628241920 & 20921948087040 & 355671353182860 & 6402052600045958 \\
\hline
$B^{9}(n)$ & 1 & 2 & 6 & 24 & 120 & 720 & 5040 & 40320 & 362880 & 3628798 & 39916758 & 479000856 & 6227008008 & 87178068432 & 1307670371280 & 20922715457280 & 355685984559360 & 6402344514209280 \\
\hline
$B^{10}(n)$ & 1 & 2 & 6 & 24 & 120 & 720 & 5040 & 40320 & 362880 & 3628800 & 39916798 & 479001554 & 6227019916 & 87178274820 & 1307674062240 & 20922784034880 & 355687312256640 & 6402371326565760 \\
\hline
$B^{11}(n)$ & 1 & 2 & 6 & 24 & 120 & 720 & 5040 & 40320 & 362880 & 3628800 & 39916800 & 479001598 & 6227020750 & 87178290164 & 1307674347420 & 20922789478080 & 355687419756480 & 6402373530940800 \\
\hline
$B^{12}(n)$ & 1 & 2 & 6 & 24 & 120 & 720 & 5040 & 40320 & 362880 & 3628800 & 39916800 & 479001600 & 6227020798 & 87178291146 & 1307674366800 & 20922789862560 & 355687427557440 & 6402373694124480 \\
\hline
$B^{13}(n)$ & 1 & 2 & 6 & 24 & 120 & 720 & 5040 & 40320 & 362880 & 3628800 & 39916800 & 479001600 & 6227020800 & 87178291198 & 1307674367942 & 20922789886624 & 355687428064992 & 6402373705032768 \\
\hline
$B^{14}(n)$ & 1 & 2 & 6 & 24 & 120 & 720 & 5040 & 40320 & 362880 & 3628800 & 39916800 & 479001600 & 6227020800 & 87178291200 & 1307674367998 & 20922789887938 & 355687428094436 & 6402373705690668 \\
\hline
$B^{15}(n)$ & 1 & 2 & 6 & 24 & 120 & 720 & 5040 & 40320 & 362880 & 3628800 & 39916800 & 479001600 & 6227020800 & 87178291200 & 1307674368000 & 20922789887998 & 355687428095934 & 6402373705726236 \\
\hline
$B^{16}(n)$ & 1 & 2 & 6 & 24 & 120 & 720 & 5040 & 40320 & 362880 & 3628800 & 39916800 & 479001600 & 6227020800 & 87178291200 & 1307674368000 & 20922789888000 & 355687428095998 & 6402373705727930 \\
\hline
$B^{17}(n)$ & 1 & 2 & 6 & 24 & 120 & 720 & 5040 & 40320 & 362880 & 3628800 & 39916800 & 479001600 & 6227020800 & 87178291200 & 1307674368000 & 20922789888000 & 355687428096000 & 6402373705727998 \\
\hline
$B^{18}(n)$ & 1 & 2 & 6 & 24 & 120 & 720 & 5040 & 40320 & 362880 & 3628800 & 39916800 & 479001600 & 6227020800 & 87178291200 & 1307674368000 & 20922789888000 & 355687428096000 & 6402373705728000 \\
\hline
\end{tabular}

\vspace{2\baselineskip}
\caption{Values $A^k(n)$ for $k,n\leq 18$. Highlighted values indicate disagreements with \cite{David1966}.}\label{table:Akn}
\begin{tabular}{|r|rrrrrrrrrrrrrrrrrr|}
\hline
 $n$ & 1 & 2 & 3 & 4 & 5 & 6 & 7 & 8 & 9 & 10 & 11 & 12 & 13 & 14 & 15 & 16 & 17 & 18 \\
\hline\hline
$A^{1}(n)$ & 1 & 0 & 0 & 0 & 0 & 0 & 0 & 0 & 0 & 0 & 0 & 0 & 0 & 0 & 0 & 0 & 0 & 0 \\
\hline
$A^{2}(n)$ & 0 & 2 & 4 & 10 & 32 & 122 & 544 & 2770 & 15872 & 101042 & 707584 & 5405530 & 44736512 & 398721962 & 3807514624 & 38783024290 & 419730685952 & 4809759350882 \\
\hline
$A^{3}(n)$ & 0 & 0 & 2 & 12 & 70 & 442 & 3108 & 24216 & 208586 & 1972904 & 20373338 & 228346522 & \textbf{2763212980} & \textbf{35926266244} & 499676669254 & 7405014187564 & 116511984902094 & 1940073930857802 \\
\hline
$A^{4}(n)$ & 0 & 0 & 0 & 2 & 16 & 134 & 1164 & 10982 & 112354 & 1245676 & 14909340 & 191916532 & \textbf{2646100822} & \textbf{38932850396} & 609137502242 & 10101955358506 & 177053463254274 & 3270694371428814 \\
\hline
$A^{5}(n)$ & 0 & 0 & 0 & 0 & 2 & 20 & 198 & 2048 & 22468 & 264538 & 3340962 & 45173518 & \textbf{652209564} & \textbf{10024669626} & 163546399460 & 2823941647390 & 51468705947590 & 987671243816650 \\
\hline
$A^{6}(n)$ & 0 & 0 & 0 & 0 & 0 & 2 & 24 & 274 & 3204 & 39420 & 514296 & 7137818 & 105318770 & \textbf{1649355338} & 27356466626 & 479446719522 & 8858271760146 & 172151975433756 \\
\hline
$A^{7}(n)$ & 0 & 0 & 0 & 0 & 0 & 0 & 2 & 28 & 362 & 4720 & 64020 & 913440 & 13760472 & 219040274 & 3681354658 & 65231186514 & 1216489698082 & 23832126613268 \\
\hline
$A^{8}(n)$ & 0 & 0 & 0 & 0 & 0 & 0 & 0 & 2 & 32 & 462 & 6644 & 98472 & 1523808 & 24744720 & 422335056 & 7575963254 & 142706934722 & 2819192544786 \\
\hline
$A^{9}(n)$ & 0 & 0 & 0 & 0 & 0 & 0 & 0 & 0 & 2 & 36 & 574 & 9024 & 145080 & 2419872 & 42129360 & 767370240 & 14631376500 & 291914163322 \\
\hline
$A^{10}(n)$ & 0 & 0 & 0 & 0 & 0 & 0 & 0 & 0 & 0 & 2 & 40 & 698 & 11908 & 206388 & 3690960 & 68577600 & 1327697280 & 26812356480 \\
\hline
$A^{11}(n)$ & 0 & 0 & 0 & 0 & 0 & 0 & 0 & 0 & 0 & 0 & 2 & 44 & 834 & 15344 & 285180 & 5443200 & 107499840 & 2204375040 \\
\hline
$A^{12}(n)$ & 0 & 0 & 0 & 0 & 0 & 0 & 0 & 0 & 0 & 0 & 0 & 2 & 48 & 982 & 19380 & 384480 & 7800960 & 163183680 \\
\hline
$A^{13}(n)$ & 0 & 0 & 0 & 0 & 0 & 0 & 0 & 0 & 0 & 0 & 0 & 0 & 2 & 52 & 1142 & 24064 & 507552 & 10908288 \\
\hline
$A^{14}(n)$ & 0 & 0 & 0 & 0 & 0 & 0 & 0 & 0 & 0 & 0 & 0 & 0 & 0 & 2 & 56 & 1314 & 29444 & 657900 \\
\hline
$A^{15}(n)$ & 0 & 0 & 0 & 0 & 0 & 0 & 0 & 0 & 0 & 0 & 0 & 0 & 0 & 0 & 2 & 60 & 1498 & 35568 \\
\hline
$A^{16}(n)$ & 0 & 0 & 0 & 0 & 0 & 0 & 0 & 0 & 0 & 0 & 0 & 0 & 0 & 0 & 0 & 2 & 64 & 1694 \\
\hline
$A^{17}(n)$ & 0 & 0 & 0 & 0 & 0 & 0 & 0 & 0 & 0 & 0 & 0 & 0 & 0 & 0 & 0 & 0 & 2 & 68 \\
\hline
$A^{18}(n)$ & 0 & 0 & 0 & 0 & 0 & 0 & 0 & 0 & 0 & 0 & 0 & 0 & 0 & 0 & 0 & 0 & 0 & 2 \\
\hline
\end{tabular}
\end{tiny}
\end{center}
\end{sidewaystable}

\end{document}